\documentclass[11pt, oneside]{article}   
\usepackage{geometry}                	
\usepackage{graphicx}	
\usepackage{amssymb}
\usepackage{float}
\usepackage[cmex10]{amsmath}
\usepackage{mathrsfs,amsthm}
\usepackage{color}

\DeclareMathOperator\erf{erf}
\DeclareMathOperator\Ei{Ei}

\newcommand\rect{\rm rect}
\newcommand\bs{\boldsymbol s}
\newcommand\bx{\boldsymbol x}
\newcommand\by{\boldsymbol y}

\newcommand\In{\operatorname{inc}}
\newcommand\Sc{\operatorname{scat}}

\newcommand\br{\boldsymbol r}

\newtheorem{remark}{Remark}

\title{Fast convolution with free-space Green's functions}
\author{Felipe Vico\footnote{
Instituto de Telecomunicaciones y Aplicaciones Multimedia (ITEAM),
Universidad Polit\` ecnica
de Val\` encia, 46022 Val\` encia, Spain. {{\em email}:\ {\tt 
{felipe.vico@gmail.com}}}, {{\sf {mferrand@dcom.upv.es.}}} }, 
Leslie Greengard\footnote{Courant Institute, New York University, New
  York, NY and Simons Foundation, New York, NY. {\em email}:
  {\tt greengard@cims.nyu.edu}}, and
Miguel Ferrando\footnotemark[1]
}

\begin{document}
\maketitle
\begin{abstract}
We introduce a fast algorithm for computing volume potentials - that is,
the convolution of a translation invariant, free-space Green's function with 
a compactly supported source distribution defined on a uniform grid.
The algorithm relies on regularizing the Fourier transform of the Green's 
function by cutting off the interaction in physical space beyond the domain 
of interest. This permits the straightforward application of trapezoidal
quadrature and the standard FFT, with superalgebraic convergence for smooth
data. Moreover, the method can be interpreted as employing a Nystrom 
discretization of the corresponding integral operator, with matrix entries which 
can be obtained explicitly and rapidly. This is of use in the design of 
preconditioners or fast direct solvers for a variety of volume integral 
equations.  The method proposed permits the computation of any derivative of 
the potential, at the cost of an additional FFT. 
\end{abstract}

\section{Introduction}

Many problems in scientific computing require the solution of a
constant coefficient elliptic partial differential equation 
subject to suitable boundary or radiation conditions.
In many cases, the free-space Green's function for the corresponding equation 
is known but involves nonlocal (long-range) interactions.
A typical example is the Helmholtz equation in ${\bf R}^d$
\[ \Delta \phi  + k^2 \phi = f, \]
where $\phi$ can be thought of as an acoustic potential and $f$ 
a known distribution of acoustic sources, which we assume to be 
supported in the bounded domain $D = [-\frac{1}{2},\frac{1}{2}]^d$. 
This can be done without loss of generality by rescaling the Helmholtz
parameter $k$. The solution which satisfies the
Sommerfeld radiation condition is well-known to be
\begin{equation}
 \phi(\bx) =
 \int_{D} g_k(\bx-\by) \, f(\by)  \, d\by.
     \label{eq:volint}
\end{equation}
where $g_k(\br) = \frac{1}{4 i} H_0(k r)$ for $d = 2$ and 
$g_k(\br) = \frac{1}{4\pi} \frac{e^{ikr}}{r}$ for $d = 3$.
Here, $r = \| \br \|_2$ and $H_0$ denotes the zeroth order Hankel 
function of the first kind.

Note that the interaction kernel is long-range, requiring 
fast algorithms to be practical, and singular at 
$r=0$, requiring accurate quadrature techniques.
In some applications, a third difficulty is that the source density $f$
is highly inhomogeneous, requiring adaptive discretization. In such
settings, intrinsically adaptive, hierarchical methods 
are required \cite{Ethridge2001,Langston2008,Biros2,mcq}.
When the density is smooth, however, and well-resolved by a uniform mesh,
it is more convenient (and generally more efficient) to use Fourier methods.  
We restrict our attention to the latter case in the present paper. 

There are two distinct ways in which Fourier methods can be applied to the 
computation of (\ref{eq:volint}). The first is the direct
discretization of the equation 
with a {\em locally-corrected} trapezoidal rule. In the two-dimensional
setting, for example, one can discretize $D$ with a uniform mesh of $N^2$ 
points and use the approximation
\begin{align*}
 \phi(nh,mh) \approx &
 \sum_{\substack{n',m' \in [-\frac{N}{2},\frac{N}{2}] \\
|n-n'|, |m-m'| > k}}
g_k((n-n')h,(m-m')h) \, f(n'h,m'h) \, h^2 \quad +  \\
 & \sum_{\substack{n',m' \in [-\frac{N}{2},\frac{N}{2}] \\
|n-n'|, |m-m'| \leq k}}
w_{n-n',m-m'} f(n'h,m'h) \ ,
\end{align*}
where $h = \frac{1}{N}$.
Several groups have shown that $k$th order accuracy can be 
achieved in this manner by a suitable choice of weights $w_{i,j}$
(see, for example, 
\cite{Aguilar2002, Aguilar2005, Singular_regular_1, Duan2009, 
goodman_1990, lowengrub_1993}).
The net sum takes the form of a discrete
(aperiodic) convolution and, hence, can be computed using the FFT with 
zero-padding in $O(N^2 \log N)$ operations. 

Alternatively, using the convolution theorem, one can write
\begin{equation}
 \phi(\bx) = \mathscr{F}^{-1}\left( \frac{F(s)}{|\bs|^2-k^2} \right)
  = \left( \frac{1}{2 \pi} \right)^d \, 
 \int_{{\bf R}^d} e^{i \bs \cdot \bx} \, \frac{F(\bs)}{|\bs|^2 - k^2} \, 
 d\bs  \, ,
     \label{eq:ift}
\end{equation}
where
\begin{equation}
F(\bs)=\mathscr{F}(f)(\bs) =
 \int_{D} e^{-i \bs \cdot \bx} f(\bx)  \, d\bx  \, .
     \label{eq:ft}
\end{equation}

$\mathscr{F}$ here denotes the Fourier transform.
The fact that $f(\bx)$ is smooth permits us to 
compute the Fourier integral in (\ref{eq:ft}) with ``spectral" accuracy.
It also ensures that the error in
truncating the Fourier integral in the inverse transform (\ref{eq:ift}) 
is rapidly decaying with $|\bs|$.
The principal difficulty in employing Fourier methods is the 
singularity $\frac{1}{|\bs|^2 - k^2}$ in the integrand.
In the case of the 
Poisson equation, this is simply $\frac{1}{|\bs|^2}$. 

It is possible to design high order rules for the inverse
Fourier transform. In the case of the
Poisson equation in three dimensions, for example, 
switching to spherical coordinates cancels the singularity entirely. 
Combining this with the nonuniform FFT yields more or less optimal schemes 
in terms of CPU time (see \cite{jiangbao} and the references therein).
This approach becomes technically more complicated for 
the Helmholtz equation, where the singularity lives on the sphere
$|\bs| = k$.

It turns out that there is a simple method that works for all long-range
Green's functions, independent of dimension, requires only the trapezoidal 
rule, achieves spectral accuracy, and is accelerated 
by the standard FFT.  Moreover, the matrix entries corresponding 
to this high order method are easily computed - a useful feature 
for either preconditioning strategies or direct solvers when using 
volume integral methods to solve variable coefficient partial 
differential equations.

Let us suppose, for the sake of simplicity,
that we seek the restriction
of the solution $\phi(\bx)$ to the unit box $D \subset {\bf R}^d$. 
Then, the maximum distance
between any source and target point in $D$ is $\sqrt{d}$.
We define 
\begin{equation}
\arraycolsep=1.4pt\def\arraystretch{2.2}
 g^L_k(\br) = \left\{
\begin{array}{cc} 
\dfrac{1}{4 i} H_0(k r) \, 
\rect\Big(\dfrac{r}{2L}\Big) & {\rm if}\ d=2 \\
\dfrac{1}{4\pi} \, \dfrac{e^{ikr}}{r}  \,
\rect\Big(\dfrac{r}{2L}\Big) & {\rm if}\ d=3 
\end{array}
\right.
\end{equation}
with 
$\rect(x)$ defined to be the characteristic function for the unit interval:
\[
\rect(x) =  \left\{
\begin{array}{cc}
 1 & {\rm for}\ |x|<1/2 \\
 0 & {\rm for}\ |x|>1/2.
\end{array}
\right.
\]
If we set $L> \sqrt{d}$ in $d$ dimensions, then 
the solution (\ref{eq:volint}) is clearly indistinguishable from
\begin{equation}
 \phi(\bx) =
 \int_{D} g^L_k(\bx-\by) \, f(\by)  \, d\by.
     \label{eq:volint_trunc}
\end{equation}

Since $g^L_k$ is compactly supported, the Paley-Wiener theorem 
implies that its Fourier transform $G^L_k$ is entire (and $C^\infty$). 
Moreover, as we shall see below, it is straightforward to compute.
In the case of the Laplace operator in three dimensions, for example,
$G^L_0 = 2 ( \frac{\sin(Ls/2)}{s})^2$. Thus, the Poisson equation
in three dimensions has the solution
\begin{equation}
 \phi(\bx) = \frac{2}{(2 \pi)^3}
 \int_{{\bf R}^3} e^{i \bs \cdot \bx} \, 
\left( \frac{\sin(L|\bs|/2)}{|\bs|} \right)^2 
F(\bs)  \, d\bs  \, .
     \label{eq:ifft}
\end{equation}
Discretization by the trapezoidal rule on the domain 
$[-\frac{N}{2},\frac{N}{2}]^d$ permits rapid evaluation using
nothing more than the FFT. The achieved accuracy is controlled by the 
rate of decay of $F(\bs)$, with spectral accuracy achieved for sufficiently
smooth $f(\bx)$ \cite{Trefethen}.

\begin{remark}
The approach described here is both elementary and quite general, 
but seems to have been overlooked in the numerical analysis literature.
An exception is the paper \cite{Vainikko2000} by Vainikko, who used volume 
Helmholtz potentials for the iterative solution of the Lippmann-Schwinger 
equation.
\end{remark}

\section{Computing the Fourier transform of truncated translation-invariant 
kernels}

Suppose now that $f(\bx)$ is a radially symmetric function:
$f(\bx)=f(r)$, where $r=|\bx|$. Then its Fourier transform $F(\bs)$ is
also radially symmetric. For $d=3$, it is easy to verify that
\begin{equation}
F(\bs)=4\pi\int_0^{\infty}\frac{\sin(sr)}{sr} \, f(r) \, r^2 \, dr
\end{equation}
where $s=|\bs|$. 
For $d=2$,
\begin{equation}
F(\bs)=2\pi\int_0^{\infty} J_0(sr) \, f(r) \, r \, dr.
\end{equation}

For the Laplace equation in three dimensions, with Green's function
\begin{equation}
\begin{aligned}
g_0^L(\bx)=\frac{1}{4\pi|\bx|} \rect\Big(\frac{|\bx|}{2L}\Big),
\end{aligned}
\end{equation}
we have
\begin{equation}
G_0^L(\bs):=\mathscr{F}\Big(g_0^L(\bx)\Big)(\bs)=
4\pi\int_0^{L}\frac{\sin(sr)}{sr} \, \frac{1}{4\pi r} \, r^2 \, dr =
2\Big(\frac{\sin(Ls/2)}{s}\Big)^2 \, ,
\end{equation}
an analytic function expressible as a power series in $s^2$.

In $\mathbb{R}^2$, where the Green's function for the Laplace equation is
\begin{equation}
g_0(\bx)=\frac{-1}{2\pi}\log |\bx|,
\end{equation}
we obtain the Fourier transform:
\begin{equation}
G_0^L(\bs):=2\pi\int_0^{+\infty}J_0(sr)g_0^L(r)rdr=
\frac{1-J_0(Ls)}{s^2}-\frac{L\log(L)J_1(Ls)}{s} \, .
\end{equation}
We set $L = 1.8 >\sqrt{3}$ in three dimensions and 
$L = 1.5 >\sqrt{2}$ in two dimensions. 
Note that, in the inverse Fourier transform 
(\ref{eq:ifft}), the frequency content of the integrand 
in the variable of integration $\bs$ is determined by 
the maximum magnitude of $\bx$, the magnitude of $L$ and $F(\bs)$
itself. It is straightforward to check that the integrand
is sufficiently sampled with a mesh that is four times finer than
in the original box: a factor of two from the fact that we are carrying
out an aperiodic convolution so that the frequency content of 
$e^{i \bs \cdot \bx} F(\bs)$ is twice greater and a factor of two from
the oscillatory behavior of the Fourier transform of the truncated kernel.
Thus, if 
the unit box is discretized using $N$ points in each dimension, we now
require a grid of size $4N$ in each dimension.
We will see in section \ref{sec:iteration} that, after a precomputation
step, this can be reduced to a factor of $2N$.

\section{Truncated kernels of mathematical physics}

We now apply the technique described above to a collection of 
Green's functions that arise in mathematical physics.
The resulting kernels in Fourier space are always $C^\infty$, as noted
above, by the Paley-Wiener theorem \cite{Hormander}.
The method is easily extended to the calculation of any derivative using 
spectral differentiation.

Tables \ref{table1} and \ref{table2} summarize the results for
various PDEs in three and two dimensions, respectively. We omit the 
derivations which are straightforward.

\begin{table}[ht]
\caption{Spectral representations of Green's Functions in 3D}
\centering
\arraycolsep=1.4pt\def\arraystretch{2.2}
\begin{tabular}{lll}  \\ [-0.3cm]
\hline
Diff. Operator & Green's function & Truncated Spectral Representation \\
\hline \\ [-0.3cm]
$\Delta$ & $g_0(r)=\dfrac{1}{4\pi r}$  &  $G^L_0(s)=2\left(\dfrac{\sin(Ls/2)}{s}\right)^2$ \\ [0.1cm]
$\Delta +k^2$ & $g_k(r)=\dfrac{e^{ikr}}{4\pi r}$  &  $G^L_k(s)=\dfrac{-1+e^{iLk}(\cos(Ls)-i\dfrac{k}{s}\sin(Ls))}{(k-s)(k+s)}$ \\ [0.1cm]
$\Delta^2$ & $g_b(r)=\dfrac{r}{8\pi}$  &  $G^L(s)=\dfrac{(2-L^2s^2)\cos(Ls)+2Ls\sin(Ls)-2}{2s^4}$ \\ [0.1cm]
$\Delta (\Delta + k^2)$ & $g_{0k}(r)=\dfrac{e^{ikr}}{4\pi r}-\dfrac{1}{4\pi r}$  &  $G^L_{0k}(s)=G^L_k(s)-G^L_0(s)$ \\ [0.1cm]
$(\Delta +\mathbf{h}\cdot\nabla)$ & $g_{\mathbf{h}}(\bx)=\dfrac{e^{i|\mathbf{h}||\bx|}}{4\pi |\bx|}e^{i\mathbf{h}\cdot\bx}$  &  $G_{\mathbf{h}}^L(\mathbf{s})=
G^L_{|\mathbf{h}|}(|\mathbf{s-h}|)$ \\ [0.1cm]
\hline
\end{tabular}
\label{table1}
\end{table}

\begin{table}[ht]
\caption{Spectral representations of Green's Functions in 2D}
\centering
{\footnotesize
\arraycolsep=1.4pt\def\arraystretch{2.2}
\begin{tabular}{lll}  \\ [-0.3cm]
\hline
Diff. Operator & Green's function &  Truncated Spectral Representation\\
\hline \\ [-0.3cm]
$\Delta$ & $g_0(r)=\dfrac{-1}{2\pi}\log (r)$  &  $G^L_0(s)=\dfrac{1-J_0(Ls)}{s^2}-\dfrac{L \log(L)J_1(Ls)}{s}$ \\ [0.1cm]
$\Delta +k^2$ & $g_k(r)=\dfrac{i}{4}H_0^{(1)}(kr)$  &  $G^L_k(s)=\dfrac{1+\dfrac{i\pi}{2}LsJ_1(Ls)H_0^{(1)}(Lk)}{s^2-k^2}$ \\ [0.1cm]
$$ & $$ & \ \ \ \ \ \ \ \ \ $-\dfrac{\dfrac{i \pi}{2}L kJ_0(Ls)H_1^{(1)}(Lk)}{s^2-k^2}$ \\ [0.1cm]
$\Delta^2$ & $g_b(r)=-\dfrac{r^2}{8\pi}\big(\log(r)-1\big) $ & $ G^L(s)=\dfrac{J_0(Ls)-1}{s^4}-\dfrac{L^3(\log(L)-1)J_1(Ls)}{4s}$ \\ [0.1cm]
$$ & $$ & \ \ \ \ \ \ \ \ \ $+\dfrac{(L\log L)J_1(Ls)}{s^3}-\dfrac{L^2(2\log L-1)J_0(Ls)}{4s^2}$ \\ [0.1cm]
$\Delta (\Delta + k^2)$ & $g_{0k}(r)=g_k(r)+g_0(r)$  &  $G^L_{0k}(s)=G^L_k(s)-G^L_0(s)$ \\ [0.1cm]
$(\Delta +\mathbf{h}\cdot\nabla)$ & $g_{\mathbf{h}}(r)=\dfrac{i}{4}H_0^{(1)}(|\mathbf{h}||\bx|)e^{i\mathbf{h}\cdot\bx}$  &  $G^L_{\mathbf{h}}(\mathbf{s})=
G^L_{|\mathbf{h}|}(|\mathbf{s-h}|)$ \\ [0.1cm]
\hline
\end{tabular}
}
\label{table2}
\end{table}

\begin{figure}[htbp]
\begin{center}
\includegraphics[width=5.9in]{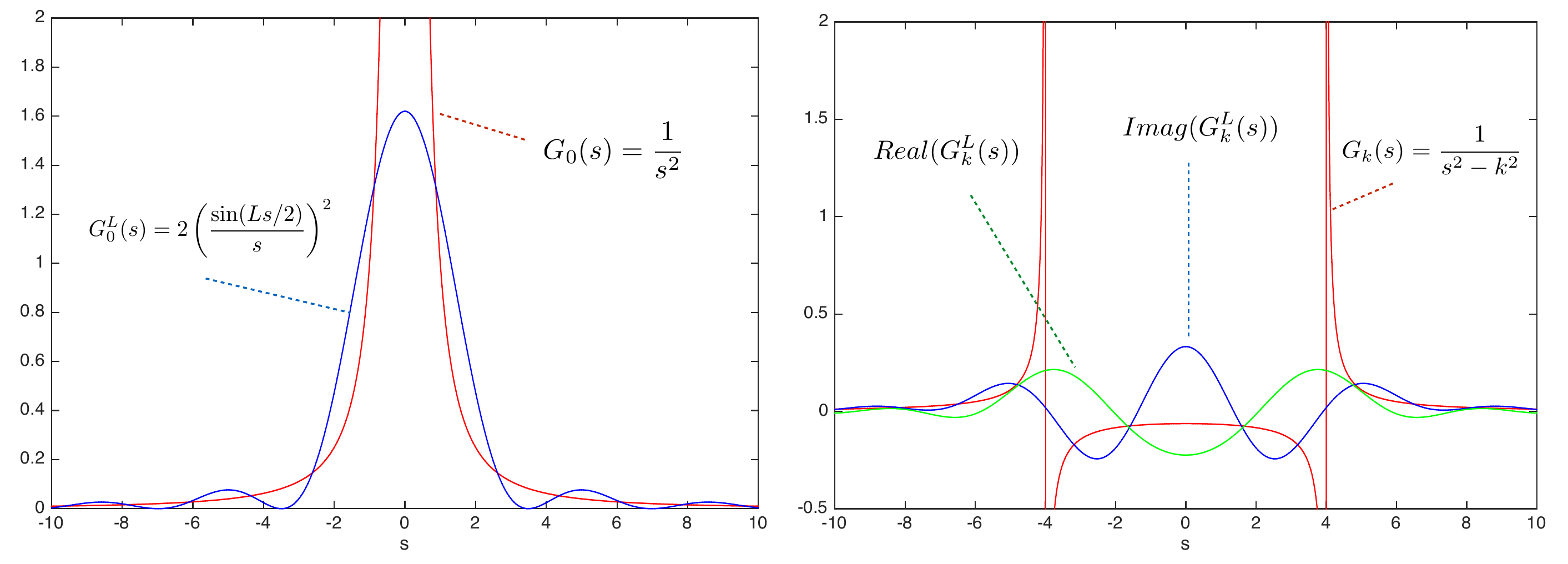} \ 
\end{center}
\caption{(l) Spectrum of the free-space Laplace kernel and the 
truncated Laplace kernel. (r) Spectrum of the free-space Helmholtz
kernel for $k=4$ and the truncated kernel with $L=1.8$. Note that 
the truncated kernels are smooth but have introduced a slight oscillation.}
\label{Lapl_spec}
\end{figure}

For illustration, we plot the spectral representations of the 
free-space and truncated Laplace and Helmholtz Green's functions
in Fig. \ref{Lapl_spec}.

\section{An explicit construction of the discretized volume integral
operators} \label{sec:iteration}

The method described above requires a grid of dimension $(4N)^d$ points
in order to compute an accurate volume integral without aliasing error.
We show now that, after a precomputation step, only an FFT of dimension 
$(2N)^d$ is required. To see this, let us consider the three dimensional
setting, with the data in the unit box denoted by 
$\rho_{ijk}=\rho(ih,jh,kh)$ where $h=1/N$ and 
$i,j,k\in\{-N/2+1,..,N/2 \}$. The solution must then take the form 
of a discrete convolution operator:
\begin{equation}
\phi_{i'j'k'} =  \sum_{i,j,k} T(i'-i,j'-j,k'-k) \rho_{ijk} \, .
\end{equation}
Thus, all entries of $T$ can be determined by simply applying the operator to 
the special right-hand side 
$\rho_{ijk}= \delta_{i0} \delta_{j0} \delta_{k0}$.
Subsequent applications of $T$ to a vector can then be carried out using standard
aperiodic convolution, which only requires zero-padding by a factor of 2.

\begin{remark}
A side effect of this precomputation is that we have generated a discrete
matrix corresponding to a high order accurate Nystr\"{o}m 
discretization of the original
volume integral operator. This is useful when implementing linear
algebraic tools such as hierarchical direct solvers, incomplete LU 
preconditioners, etc.
\end{remark}

\begin{figure}[htbp]
\begin{center}
\includegraphics[width=4in]{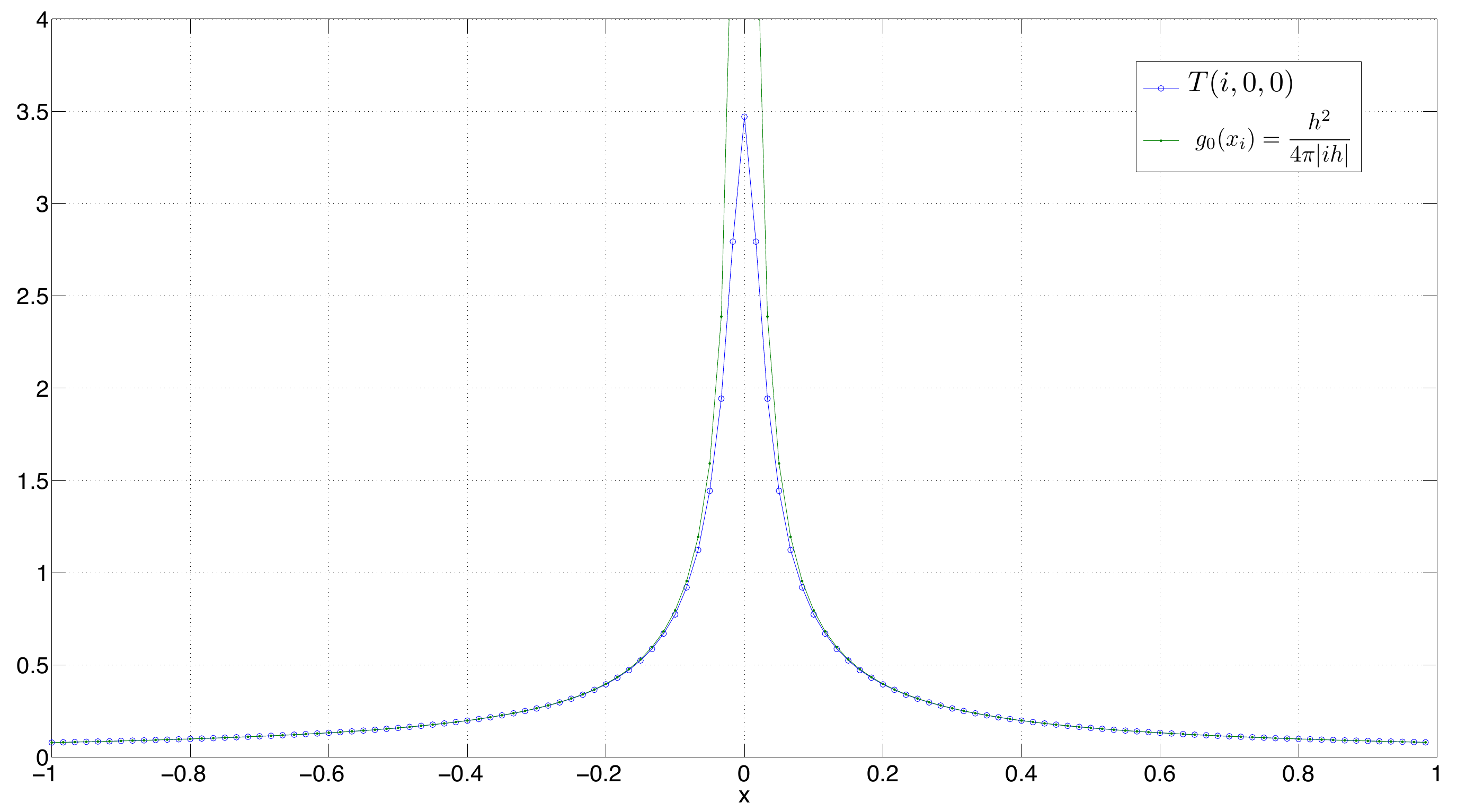}%
\end{center}
\caption{Comparing the naive trapezoidal rule on the original Green's
function and the high order mollified Green's function along the 
line $y=z=0$.}
\label{Coef_1}
\end{figure}

It is worth plotting the resulting entries of $T$ and comparing them to a
naive trapezoidal approximation (which blows up when $i=i'$, $j=j'$ and
$k=k'$). As can be seen in Fig. \ref{Coef_1}, our high order rule takes the 
form of a mollified Green's function - with no significant oscillations in 
sign or other difficulties that plague many high order quadrature generation
techniques.

\section{Numerical results}

In this section, we illustrate the performance of the method described above.
Our first examples simply involve convolution of the free-space 
Green's function
with a Gaussian source, since the exact solution is available analytically.
We also solve a variable dielectric Poisson-Boltzmann equation
and a Lippmann-Schwinger type integral equation for variable medium 
scattering problems.

\subsection{Convolution with a Gaussian source}

Suppose now that, in three dimensions,
the source distribution is given by a simple Gaussian:
\begin{equation}
\rho(r):=\frac{1}{\sigma^3(2\pi)^{3/2}}e^{-\frac{r^2}{2\sigma^2}} \, .
\end{equation}
Then, the solution to the Poisson equation is given by
\begin{equation}\label{test_lapl_3D}
[g_0\ast\rho](\bx)=\frac{1}{4\pi r}\erf\Big(\frac{r}{\sqrt{2}\sigma}\Big) \, .
\end{equation}
For the Helmholtz equation, the solution is a little more complicated
but also straightforward to compute:
\begin{equation}\label{test_helm_3D}
[g_k\ast\rho](\bx)=\frac{1}{4\pi r}e^{-\frac{\sigma^2k^2}{2}}
\left[ Real\left(e^{-ikr}\erf\left(\frac{2\sigma^2ik-2r}{2\sqrt{2\sigma^2}}\right) \right)-i\sin(kr)   \right] \, .
\end{equation}
For the biharmonic equation, we have
\begin{equation}\label{test_biharm_3D}
[g_b\ast\rho ](\bx)= \frac{1}{8\pi} \left[ \sigma \sqrt{\frac{2}{\pi}} +
\erf{\Big(\frac{r}{\sigma\sqrt{2}}\Big)}\Big(\frac{\sigma^2}{r}+r\Big) \right].
\end{equation}
Similarly, in two dimensions, with
\begin{equation}
\rho(r):=\frac{1}{2\pi\sigma^2}e^{-\frac{r^2}{2\sigma^2}} \, ,
\end{equation}
we have the following solution for the Poisson equation:
\begin{equation}\label{test_lapl_2D}
[g_0\ast\rho](\bx)=\frac{-1}{4\pi}\left[
\Ei\Big(\frac{r^2}{2\sigma^2}\Big)+\log(r^2)\right].
\end{equation}
For the Helmholtz equation, we have
\begin{equation}\label{test_helm_2D}
[g_k\ast\rho](\bx)=
\frac{H_0(kr)}{4\sigma^2}
\int_0^r J_0(ky) \, e^{-\frac{y^2}{2\sigma^2}} \, y dy+
\frac{J_0(kr)}{4\sigma^2}
\int_r^{+\infty}H_0(kr) \, e^{-\frac{y^2}{2\sigma^2}} \, y dy 
\end{equation}
and for the biharmonic equation, we have
\begin{equation}\label{test_biharm_2D}
[g_b\ast\rho ] (\bx)=-\frac{\sigma^2}{8\pi}
\left[ \left(\frac{r^2}{2\sigma^2}+1 \right)
\widetilde{\Ei}\left(\frac{r^2}{2\sigma^2} \right)-e^{-\frac{r^2}{2\sigma^2}} 
\right]+c_2r^2+c_1 \, r,
\end{equation}
where
\begin{equation}
\begin{aligned}
\widetilde{\Ei}(x):=&\Ei(x)+\log(x)+\gamma\\
c_1:=&\frac{\sigma^2}{8\pi}\Big(\gamma+\log\big(\frac{1}{2\sigma^2}\big)\Big)\\
c_2:=&\frac{1}{8\pi}\Big(\frac{\gamma}{2}+\frac{1}{2}\log\big(\frac{1}{2\sigma^2}\big)+1\Big) \, .
\end{aligned}
\end{equation}

\begin{figure}[htbp]
\begin{center}
\includegraphics[width=5.5in]{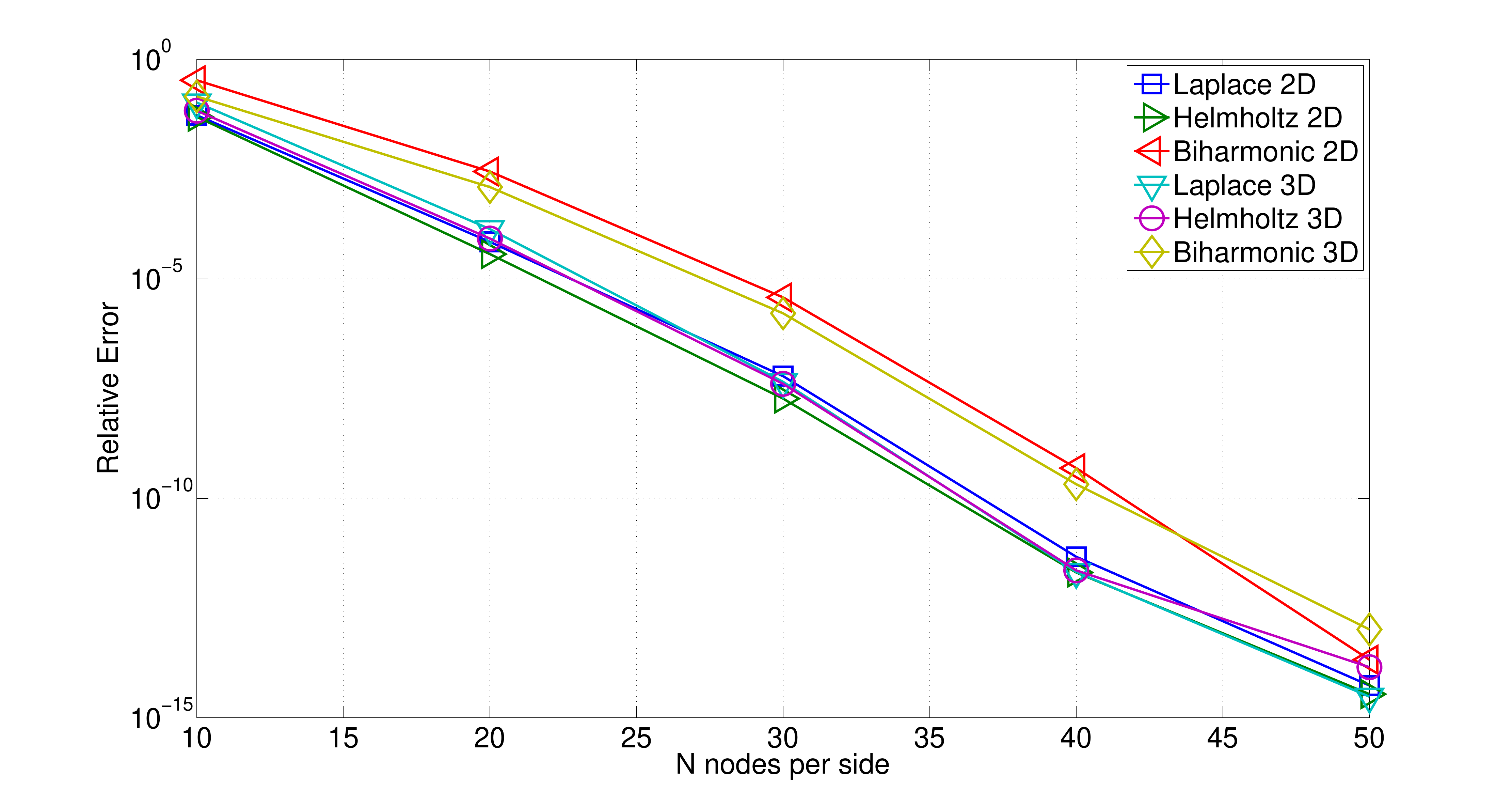}%
\end{center}
\caption{Convergence of the truncated Green's function Fourier method
in solving the Poisson, Helmholtz and biharmonic equations in two and three
dimensions for a single Gaussian source with $\sigma=0.05$
(see eqs. (\ref{test_lapl_3D}-\ref{test_biharm_3D}),
(\ref{test_lapl_2D}-\ref{test_biharm_2D})). 
The Helmholtz parameter was set to $k=2$.}
\label{Error_Convolution}
\end{figure}

In Fig. \ref{Error_Convolution}, we plot the error in each of these solutions
when computed using the truncated Green's function Fourier method. 
Spectral accuracy is evident in each case.

\subsection{Non-oscillatory elliptic equations with variable coefficients}

A variety of problems in computational physics require the solution of the
divergence-form elliptic partial differential equation
\begin{equation}
\nabla\cdot\epsilon(\bx)\nabla\phi - \lambda^2 \phi =\rho(\bx)
\end{equation}
where $\epsilon$ is a known, smooth perturbation of a background constant
$\epsilon+0$, and where both $\rho(\bx)$ and $\epsilon-\epsilon_0$ have 
compact support. 

\begin{figure}[htbp]
\begin{center}
\includegraphics[width=5in]{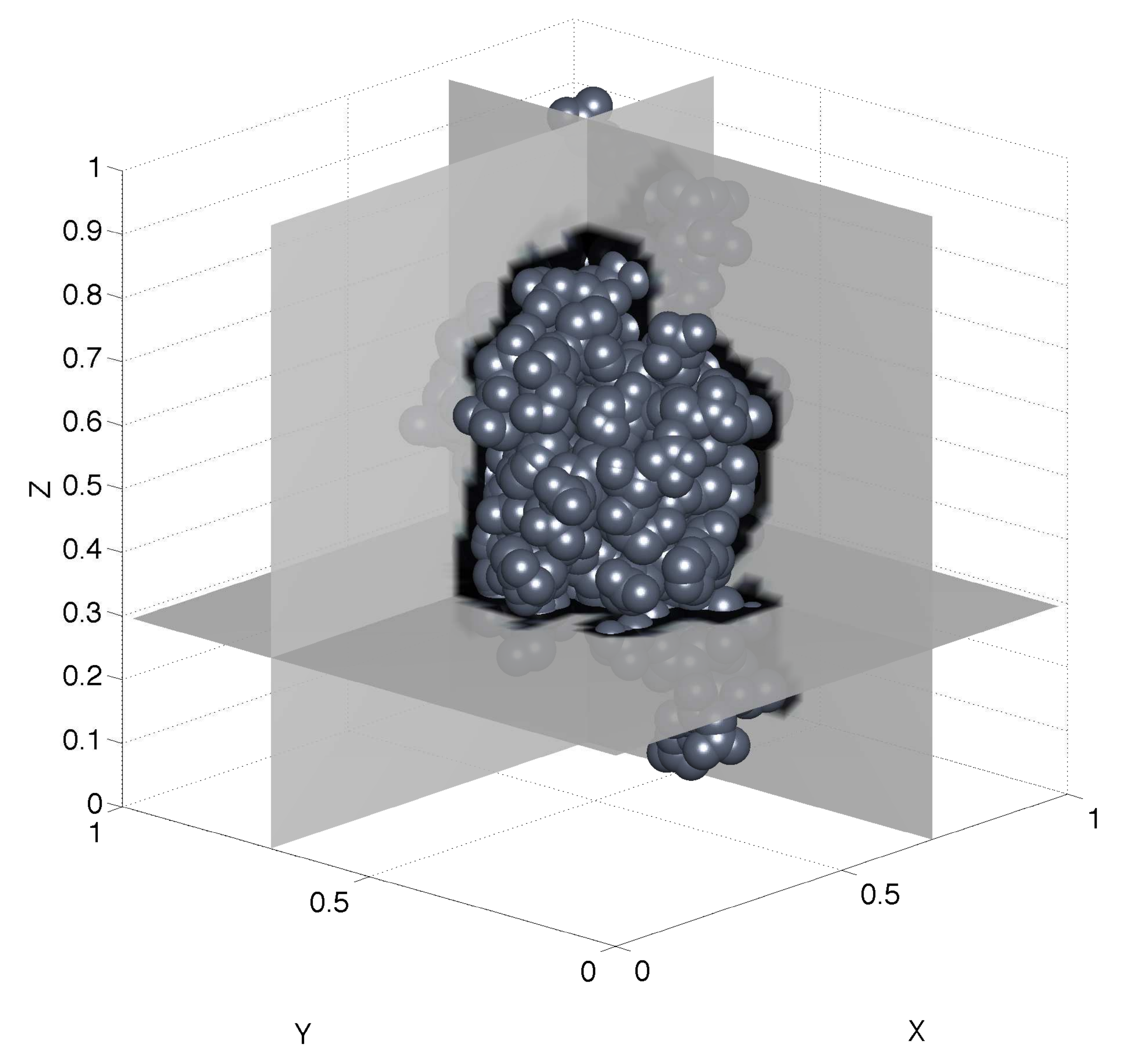}%
\end{center}
\caption{A protein molecule and the associated smooth dielectric function.}
\label{Protein_3D}
\end{figure}

In molecular electrostatics, this equation is referred to as the linearized
Poisson-Boltzmann equation. While most models make use of a sharp dielectric
interface (so that $\epsilon$ is piecewise constant), there is also 
interest in using approximations of the dielectric that vary smoothly
\cite{DelPhi,Nicholls}. Following these references, 
a formula for $\epsilon$ is determined by first
assuming we are given a macromolecule with 
$M$ atoms, with a density
\[  \alpha_i(\bx) = exp[ -r_i^2/(\mu^2R_i^2)] 
\]
centered on the $i$th atom, where $r_i$ denotes the distance of $\bx$ from
the atomic center,
$R_i$ is the van der Waals radius of the atom and $\mu$ is a user-specified
variance. From this, a function 
\[ q(\bx) = 1 - \prod_i [1 - \alpha_i(\bx)]
\] 
is constructed and, finally, 
\[ \epsilon(\bx) = q(\bx) \, \epsilon_{in} + (1- q(\bx) \,  \epsilon_{out}. 
\] 
In our example, we let $M = 1235$ with, $R_i = 0.022$ and 
$\mu^2 = 2$.
Figure \ref{Protein_3D} shows the molecule represented as a union of
spheres and the associated dielectric function $\epsilon(\bx)$ in darker gray.

Restricting our attention to the case $\lambda = 0$ for the sake of simplicity,
we may represent the solution in the form
\begin{equation}
\phi(\bx)=\int_{\mathbb{R}^3}\frac{1}{4\pi|\bx-\by|}\sigma(\by)dV_{\by} \, .
\end{equation}
This leads directly to the following
second kind Fredholm equation for the unknown density $\sigma(\bx)$:
\begin{equation}
- \epsilon(\bx) \sigma+\nabla\epsilon(\bx) \cdot\nabla\int_{\mathbb{R}^3}\frac{1}{|\bx-\by|}\sigma(\by)dV_{\by}=\rho(\bx) \, .
\label{pbeq}
\end{equation}

We discretize $\sigma$ on a uniform mesh with $N^3$ points 
and use the truncated Green's function
Fourier method described above to convert (\ref{pbeq}) into a dense system of 
equations which we solve iteratively using GMRES. Each matrix vector product
requires $O(N^3 \log N)$ operations using the FFT.

Our results are summarized in Table \ref{table_Laplace}. 
$N_{tot} = N^3$ denotes the total number of unknowns, $N_{iter}$ denotes
the total number of GMRES iterations, 
$E_2$ denotes the relative error with respect to the reference solution 
in $L_2$ for $N=250$, $E_{inf}$ denotes the relative error
in $L_\infty$, and $T_{solve}$ denotes the solution
time in seconds required on a workstation with 
two Intel Xeon E5-2450 processors with 8 cores per processor and
64 GB of memory. $T_{precomp}$ denotes the time required for 
precomputation, as discussed in section 
\ref{sec:iteration}, which requires a single FFT of dimension $(4N)^3$.

\begin{table}[ht]
\caption{Fast, iterative solution of the linearized Poisson-Boltzmann equation 
using the truncated Green's function Fourier method.}
{\small
\centering
\begin{tabular}{ccccccc}  \\ [-0.3cm]
\hline
$N_{tot}$ & $N$ & $E_2$ & $E_{inf}$ & $N_{iter}$ & $T_{solve}$ $(s)$ & $T_{precomp}$ (s)\\
\hline \\ [-0.3cm]
 $1000000$ & $100$ & $2.08\times10^{-8}$ & $2.54\times10^{-6}$ & $12$ & $14.3$ & $4.4$\\ [0.1cm]
 $3375000$ & $150$ & $2.74\times10^{-10}$ & $3.97\times10^{-8}$ & $16$ & $76.9$ & $13$ \\ [0.1cm]
 $8000000$ & $200$ & $7.92\times10^{-12}$ & $1.06\times10^{-9}$ & $16$ & $205$ &  $28.9$\\ [0.1cm]
 $15625000$ & $250$ & $-$ & $-$ & $16$ & $421$ &  $75.8$\\ [0.1cm]
\hline \\ [-0.3cm]
\end{tabular}
\label{table_Laplace}
}
\end{table}
\subsection{Lippmann-Schwinger equation for wave scattering}

In our last set of examples, we study the performance of the 
Lippmann-Schwinger integral equation for solving variable coefficient 
scattering problems
in $\mathbb{R}^2$ and $\mathbb{R}^3$. 
The governing equation is a Helmholtz equation of the form
\begin{equation}
\Delta \phi^{\Sc}+k^2(1+q(\bx))\phi^{\Sc}=-k^2q(\bx)\phi^{\In}
\end{equation}
where $\phi^{\Sc}$ is assumed to satisfy the usual Sommerfeld radiation
condition.
We assume $q(\bx)$ has compact support. Using a volume integral
representation for the solution in $\mathbb{R}^2$:
\begin{equation}
\phi^{\Sc}(\bx)=\int_{D} H_0(k |\bx-\by|) \, \sigma(\by) \, 
dV_{\by} \, ,
\end{equation}
we obtain the second kind integral equation
\begin{equation}
-\sigma+k^2q(\bx)\int_{D} H_0(k|\bx-\by|) \, \sigma(\by) \, dV_{\by}=
-k^2q(\bx)\phi^{\In} \, .
\end{equation}
Similarly, in $\mathbb{R}^3$ we get:
\begin{equation}
-\sigma+k^2q(\bx)\int_{D} \frac{e^{ik|\bx-\by|}}{4\pi |\bx-\by|} \, \sigma(\by) \, dV_{\by}=
-k^2q(\bx)\phi^{\In} \, .
\end{equation}
(This is the dual of the usual Lippmann-Schwinger equation.)

\begin{figure}[htbp]
\begin{center}
\includegraphics[width=5.9in]{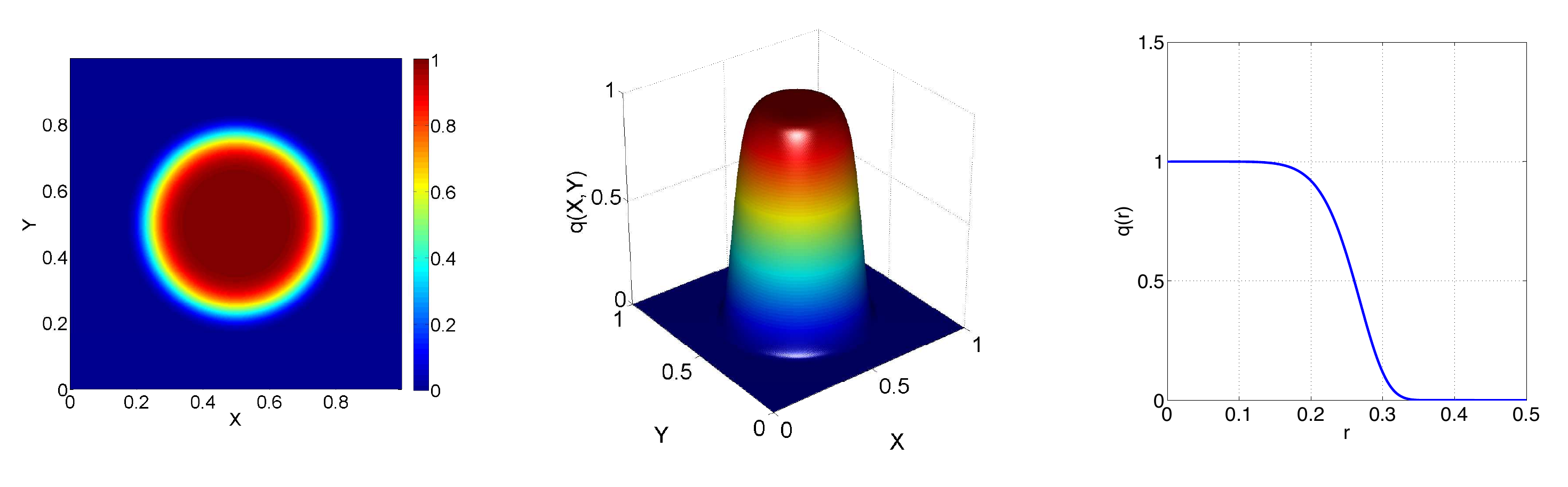}
\end{center}
\caption{Contrast function $q(\bx)$ for the smoothly filtered disk
plotted as colored contours (left), as a surface (center), and as a function
of radius (right).}
\label{SCFB_indall}
\end{figure}

\begin{figure}[htbp]
\begin{center}
\includegraphics[width=5.9in]{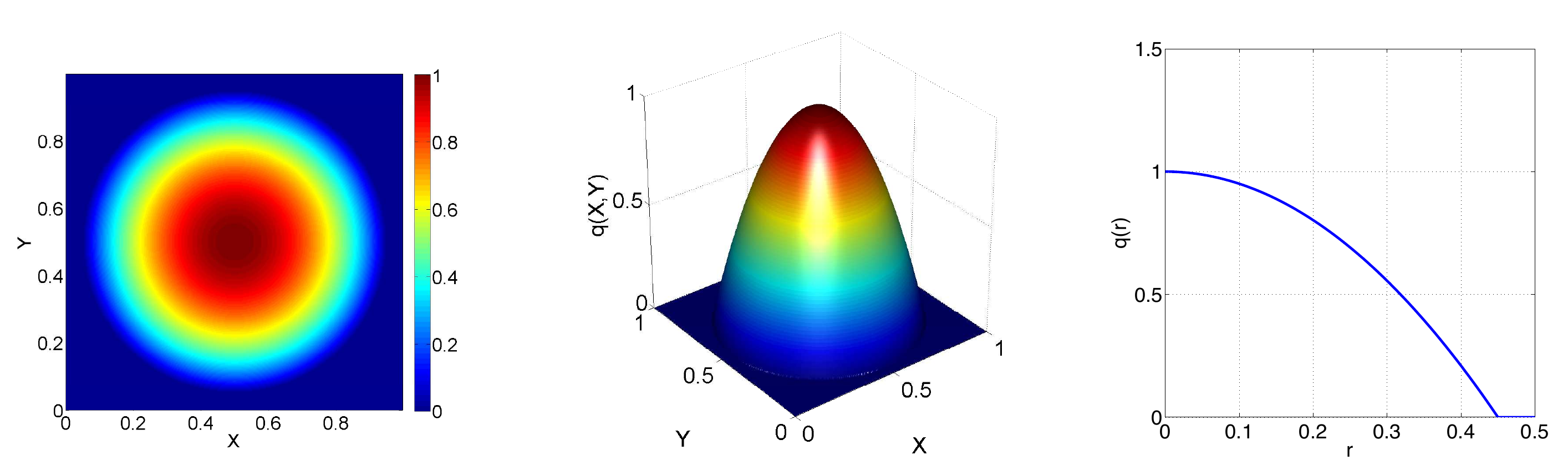}
\end{center}
\caption{Contrast function $q(\bx)$ for the Luneburg lens
plotted as colored contours (left), as a surface (center), and as a function
of radius (right).}
\label{Luneb_indexall}
\end{figure}

\begin{figure}[htbp]
\begin{center}
\includegraphics[width=5.9in]{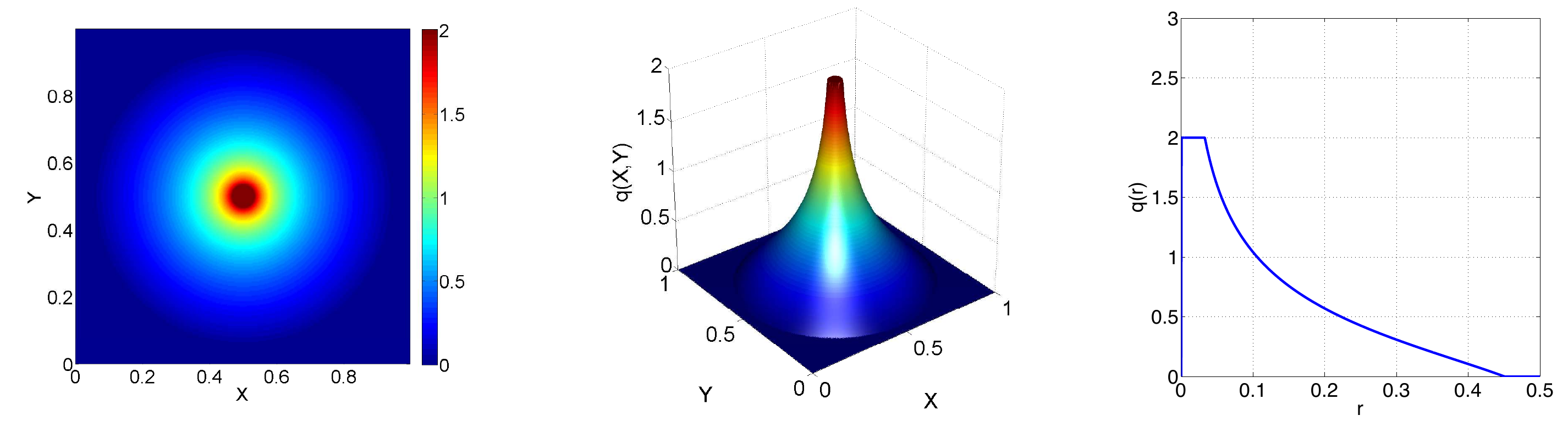}
\end{center}
\caption{Contrast function $q(\bx)$ for the Eaton lens
plotted as colored contours (left), as a surface (center), and as a function
of radius (right).}
\label{Eaton_indexall}
\end{figure}

We consider four cases: a smoothly filtered flat dielectric disk in 2D,
the 2D ``Luneburg" lens, the 2D 
``Eaton" lens and a smoothed dielectric cube in 3D.
The smoothly filtered disk 
(Fig. \ref{SCFB_indall}) has a contrast function given by
\begin{equation}
q(\bx)=e^{-\frac{1}{2}\big(\frac{|\bx|}{0.25}\big)^8}\, .
\end{equation}
The Luneburg lens 
(Fig. \ref{Luneb_indexall}) 
is designed to focus an incoming wave to a single point
\cite{luneburg1964}, with $q$ given by
\begin{equation}
q(\bx)=1-\Big(\frac{|\bx|}{0.45}\Big)^2 \, .
\end{equation}
The Eaton lens 
(Fig. \ref{Eaton_indexall}) 
is designed to bend light through an angle \cite{Eaton}, with 
$q(\bx) = n^2(\bx)-1$, where the 
refractive index $n$ is given by the implicit equation
\begin{equation}
n^2(\bx) =\frac{0.45}{n(\bx)|\bx|}+\sqrt{\Big(\frac{0.45}{n(\bx)|\bx|}\Big)^2-1}
\end{equation}
In order to avoid the blowup in $n$ at the origin, the refractive index is 
truncated at a maximum value of $n_{\max}=\sqrt{3}$, corresponding to
$q_{\max}=2$. Finally,
the smoothed cube has a contrast function given by:
\begin{equation}
q(\bx)=e^{-\frac{1}{2}\Big(\big(\frac{x}{0.25}\big)^8+\big(\frac{y}{0.25}\big)^8+\big(\frac{z}{0.25}\big)^8\Big)}\, .
\end{equation}

We solve each Lippmann-Schwinger equation iteratively, using
Bi-CGStab with a tolerance of $10^{-12}$ 
for the iteration, since it has minimal storage requirements.
This requires two matrix-vector products per iteration, each involving
two applications of the FFT using the 
truncated Green's function Fourier method.

Except for the Eaton lens, the incoming wave is chosen to be a plane
wave propagating to the right. The incoming wave for the Eaton lens is
given by a  Gaussian beam of the form
\[
\phi^{\In}= \overline{H_0(kR)}e^{-0.5k}
\]
where 
\[ R=\sqrt{(x-x_0^c)^2+(y-y_0^c)^2},\ 
x_0^c =-0.01- 0.5i, y_0^c = 0.77 \, .
\]

\begin{figure}[htbp]
\begin{center}
\includegraphics[width=3.5in]{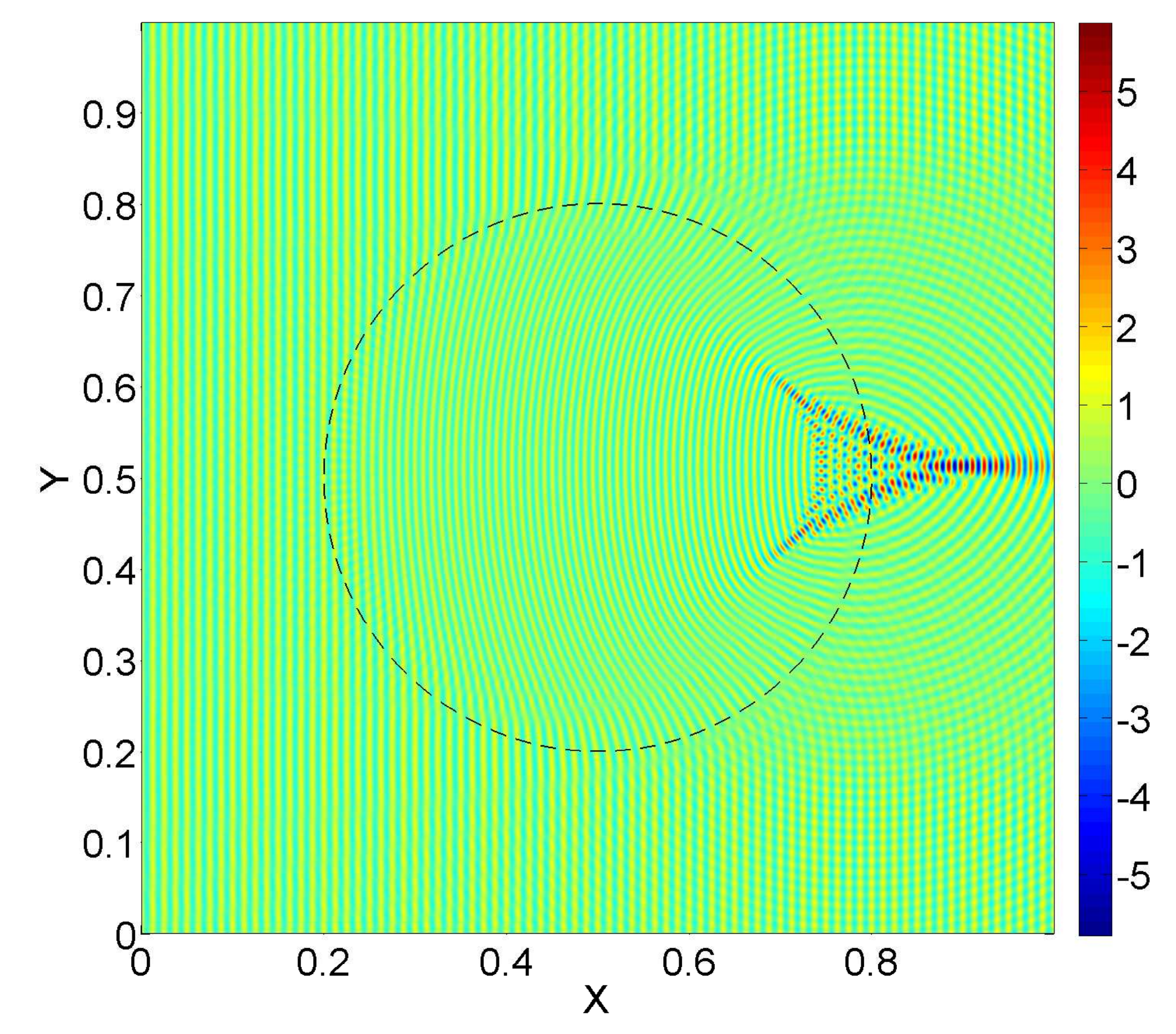}\ 
\end{center}
\caption{Scattering from a smoothly filtered disk with 
radius $R=40 \lambda_0$. 
in a unit square of size $80\lambda_0\times 80\lambda_0$.
We compute a reference solution with 
$N_{tot}=40960000 = 6400^2$ points.
\label{SCFB_field}}
\end{figure}

\begin{figure}[htbp]
\begin{center}
\includegraphics[width=5.9in]{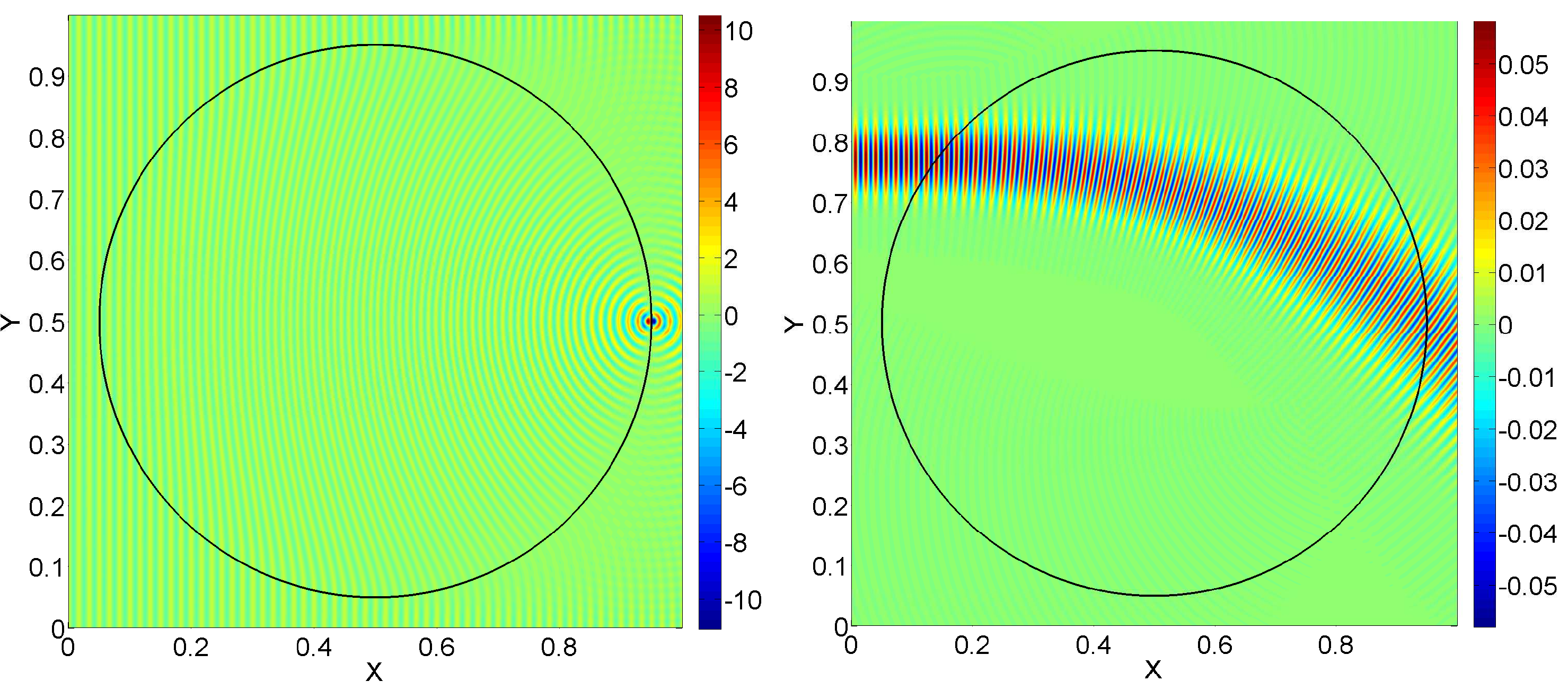}\ 
\end{center}
\caption{(left) Scattering by the two-dimensional Luneburg Lens,
(right) Bending of an incoming beam by the two-dimensional Eaton Lens. 
In both cases,
the lens radius is $R=27 \lambda_0$ in a
unit square of size $60\lambda_0\times 60\lambda_0$.
The reference solution was computed with $N_{tot}=40960000 = 6400^2$ points.
\label{Luneb_field}}
\end{figure}

\begin{figure}[H]
\begin{center}
\includegraphics[width=4.8in]{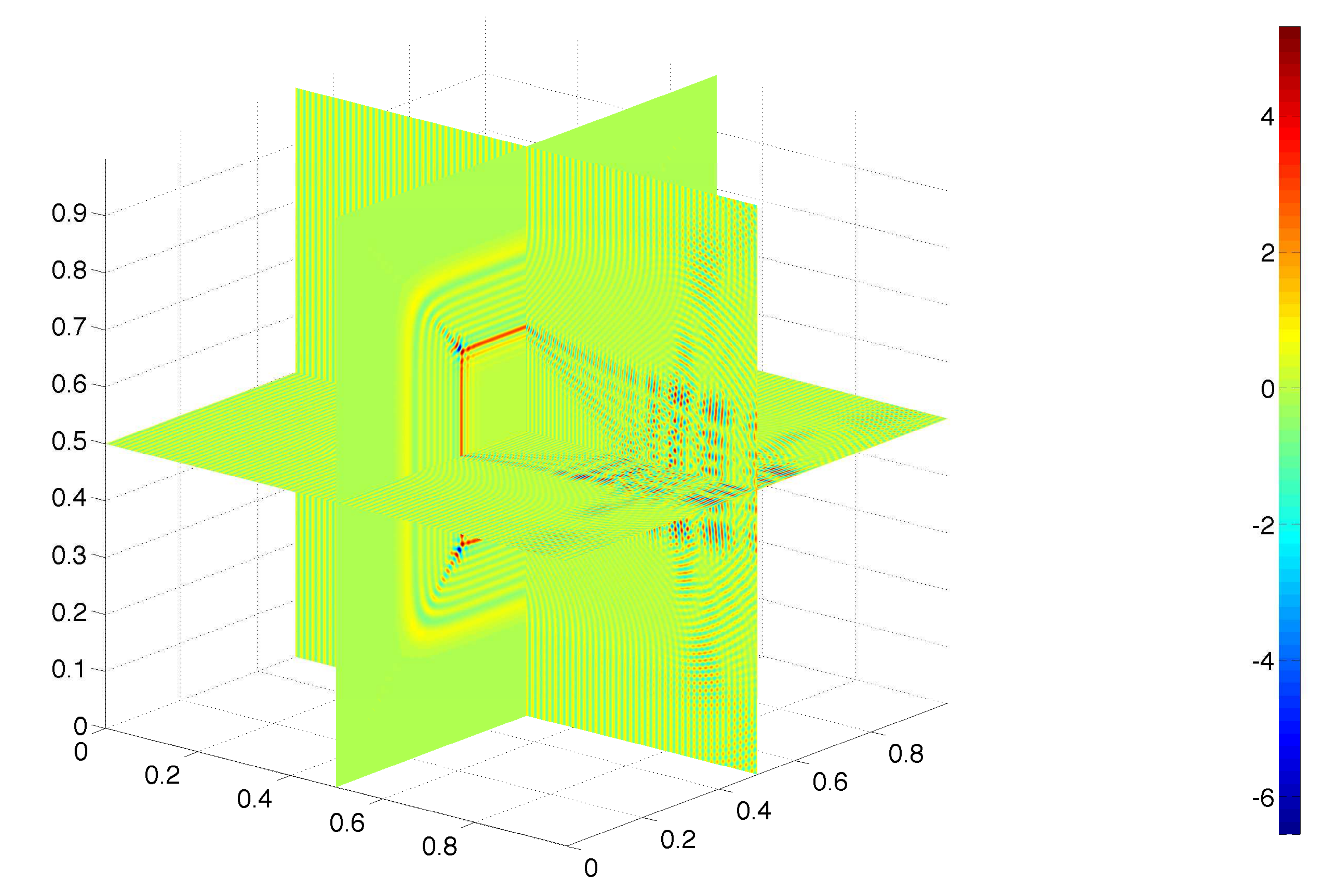}%
\end{center}
\caption{Three-dimensional scattering by a smoothed cube. 
We discretize a unit cell of dimension 
$80\lambda_0\times 80\lambda_0\times 80\lambda_0$ and
compute a reference solution with 
$N_{tot}=300^3$ points.}
\label{SSFB}
\end{figure}

Tables \ref{table_SCFB}-\ref{table_SsquareFB} show timings and errors
for various frequencies and discretizations, while Figs.
\ref{SCFB_field}-\ref{SSFB} show the computed solution.
In these tables, size denotes the dimensions of the unit box in wavelengths,
$N_{tot}$ denotes the total number of points in the discretization,
$N$ denotes the number of points in a linear dimension, 
$E_2$ denotes the relative error with respect to the reference solution 
in $L_2$, $E_{inf}$ denotes the relative error
in $L_\infty$, and $N_{matvec}$ denotes
the total number of matrix-vector products needed in the 
Bi-CGStab iteration.
As above,
$T_{solve}$ denotes the solution time in 
seconds on a workstation with 
two Intel Xeon E5-2450 processors with 8 cores per processor and
64 GB of memory, and $T_{precomp}$ denotes the time required for 
precomputation, as discussed in section 
\ref{sec:iteration}. A reference solution is computed using 
$6400 \times 6400$ points in the two-dimensional examples and using 
$300 \times 300 \times 300$ points in the three-dimensional example.

Note that spectral convergence rates are
evident for smooth dielectric contrast functions. For the 
non-smooth Eaton and Luneburg lenses, the numerical convergence rate
is closer to second order accuracy but with a small constant, so that
high precision is achieved with a modest number of points per wavelength.

\vspace{.2in}

\begin{table}[H]
\caption{Data for the smoothly filtered
disk in two dimensions (see text for discussion). Timings are in seconds.}

\vspace{.2in}

{\small
\centering
\begin{tabular}{cccccccc}  \\ [-0.3cm]
\hline
Size $(\lambda_{0})$ & $N_{tot}$ & $N$ & $E_2$ & $E_{inf}$ & $N_{matvec}$ & $T_{solve}$ & $T_{precomp}$ \\
\hline \\ [-0.3cm]
$1$ & $400$ & $20$ & $1.4\times10^{-4}$ & $2.1\times10^{-4}$ & $17$ & $0.382$ & $0.222$\\ [0.1cm]
$1$ & $2500$ & $50$ & $3.2\times10^{-8}$ & $3.2\times10^{-8}$ & $15$ & $0.388$ & $0.225$ \\ [0.1cm]
$1$ & $10000$ & $100$ & $8.7\times10^{-13}$ & $1.1\times10^{-12}$ & $15$ & $0.632$ &  $0.5$\\ [0.1cm]
$1$ & $ 40960000$ & $6400$ & $-$ & $-$ & $15$ & $149$ &  $152$\\ [0.1cm]
\hline
$20$ & $6400$ & $80$ & $4.2\times10^{-5}$ & $6.7\times10^{-5}$ & $332$ & $1.38$ & $0.26$\\ [0.1cm]
$20$ & $10000$ & $100$ & $4.5\times10^{-8}$ & $8.8\times10^{-8}$ & $333$ & $1.83$ & $0.40$ \\ [0.1cm]
$20$ & $19600$ & $140$ & $4.1\times10^{-11}$ & $6.3\times10^{-11}$ & $335$ & $2.39$ &  $0.303$\\ [0.1cm]
$20$ & $ 40960000$ & $6400$ & $-$ & $-$ & $335$ & $3170$ &  $143$\\ [0.1cm]
\hline
$80$ & $62500$ & $250$ & $6.7\times10^{-5}$ & $1.0\times10^{-4}$ & $2938$ & $58.5$ & $0.503$\\ [0.1cm]
$80$ & $72900$ & $270$ & $1.2\times10^{-7}$ & $2.2\times10^{-7}$ & $2990$ & $67.1$ & $0.518$ \\ [0.1cm]
$80$ & $102400$ & $320$ & $1.6\times10^{-10}$ & $2.8\times10^{-10}$ & $2906$ & $83.1$ &  $0.61$\\ [0.1cm]
$80$ & $ 40960000$ & $6400$ & $-$ & $-$ & $2948$ & $29283$ &  $150$\\ [0.1cm]
\hline
\end{tabular}
\label{table_SCFB}
}
\end{table}

\begin{table}[H]
\caption{Data for the two-dimensional
Luneburg lens, where $\lambda_0$ denotes
the {\em free-space} wavelength. 
A reference solution is computed using a
$6400 \times 6400$ grid (except for $\lambda_0 = 60$), where
a $3200 \times 3200$ grid is used. Timings are in seconds.}
\vspace{.2in}

{\small
\centering
\begin{tabular}{cccccccc}  \\ [-0.3cm]
\hline
Size $(\lambda_{0})$ & $N_{tot}$ & $N$ & $E_2$ & $E_{inf}$ & $N_{matvec}$ & $T_{solve}$ & $T_{precomp}$ \\
\hline \\ [-0.3cm]
$1$ & $640000$ & $800$ & $1.26\times10^{-7}$ & $2.89\times10^{-7}$ & $17$ & $ 2.41$ & $2.09$ \\ [0.1cm]
$1$ & $2560000$ & $1600$ & $2.09\times10^{-8}$ & $5.09\times10^{-8}$ & $17$ & $6.86$ &  $8.83$\\ [0.1cm]
$1$ & $10240000$ & $3200$ & $2.93\times10^{-9}$ & $7.50\times10^{-9}$ & $17$ & $40.3$ &  $ 34.9$\\ [0.1cm]
\hline
$20$ & $640000$ & $800$ & $3.18\times10^{-5}$ & $4.01\times10^{-5}$ & $582$ & $67.3$ & $2.33$ \\ [0.1cm]
$20$ & $2560000$ & $1600$ & $5.84\times10^{-6}$ & $8.42\times10^{-6}$ & $581$ & $243$ &  $9.02$\\ [0.1cm]
$20$ & $10240000$ & $3200$ & $8.45\times10^{-7}$ & $1.39\times10^{-6}$ & $590$ & $1190$ &  $34.3$\\ [0.1cm]
\hline
$40$ & $640000$ & $800$ & $7.53\times10^{-5}$ & $8.59\times10^{-5}$ & $1415$ & $163$ & $2.37$ \\ [0.1cm]
$40$ & $2560000$ & $1600$ & $1.6\times10^{-5}$ & $1.88\times10^{-5}$ & $1393$ & $740$ &  $9.24$\\ [0.1cm]
$40$ & $10240000$ & $3200$ & $3.21\times10^{-6}$ & $3.89\times10^{-6}$ & $1321$ & $2969$ &  $34.5$\\ [0.1cm]
\hline
\hline
$60$ & $640000$ & $800$ & $1.26\times10^{-4}$ & $1.40\times10^{-4}$ & $3844$ & $449$ & $2.45$ \\ [0.1cm]
$60$ & $ 2250000$ & $1500$ & $3.13\times10^{-5}$ & $3.54\times10^{-5}$ & $3482$ & $1611$ &  $8.27$\\ [0.1cm]
$60$ & $10240000$ & $3200$ & $-$ & $-$ & $5220$ & $12322$ &  $34.5$\\ \hline
\end{tabular}
\label{table_Luneburg}
}
\end{table}

\begin{table}[H]
\caption{Data for the two-dimensional
Eaton lens, where $\lambda_0$ denotes
the {\em free-space} wavelength. 
A reference solution is computed using a
$6400 \times 6400$ grid and timings are in seconds.}

\vspace{.2in}

{\small
\centering
\begin{tabular}{cccccccc}  \\ [-0.3cm]
\hline
Size $(\lambda_{0})$ & $N_{tot}$ & $N$ & $E_2$ & $E_{inf}$ & $N_{matvec}$ & $T_{solve}$ & $T_{precomp}$ \\
\hline \\ [-0.3cm]
$1$ & $640000$ & $800$ & $1.01\times10^{-7}$ & $1.62\times10^{-7}$ & $15$ & $1.82$ & $1.85$ \\ [0.1cm]
$1$ & $2560000$ & $1600$ & $9.36\times10^{-8}$ & $1.75\times10^{-7}$ & $15$ & $6.63$ &  $8.74$\\ [0.1cm]
$1$ & $10240000$ & $3200$ & $3.53\times10^{-9}$ & $6.86\times10^{-9}$ & $15$ & $35.5$ &  $34.4$\\ [0.1cm]
\hline
$20$ & $640000$ & $800$ & $4.96\times10^{-6}$ & $1.54\times10^{-5}$ & $388$ & $51.1$ & $2.49$ \\ [0.1cm]
$20$ & $2560000$ & $1600$ & $8.50\times10^{-7}$ & $2.96\times10^{-6}$ & $390$ & $203$ &  $9.10$\\ [0.1cm]
$20$ & $10240000$ & $3200$ & $1.25\times10^{-7}$ & $4.77\times10^{-7}$ & $386$ & $850$ &  $34.8$\\ [0.1cm]
\hline
$40$ & $640000$ & $800$ & $1.4\times10^{-5}$ & $3.67\times10^{-5}$ & $958$ & $105$ & $2.26$ \\ [0.1cm]
$40$ & $2560000$ & $1600$ & $2.34\times10^{-6}$ & $8.77\times10^{-6}$ & $956$ & $419$ &  $10.4$\\ [0.1cm]
$40$ & $10240000$ & $3200$ & $3.98\times10^{-7}$ & $1.55\times10^{-6}$ & $968$ & $2163$ &  $40$\\ [0.1cm]
\hline
$60$ & $640000$ & $800$ & $2.81\times10^{-5}$ & $6.90\times10^{-5}$ & $2064$ & $276$ & $2.68$ \\ [0.1cm]
$60$ & $2560000$ & $1600$ & $4.61\times10^{-6}$ & $1.45\times10^{-5}$ & $2038$ & $1065$ &  $9.06$\\ [0.1cm]
$60$ & $10240000$ & $3200$ & $7.46\times10^{-7}$ & $2.87\times10^{-6}$ & $2024$ & $4550$ &  $33$\\ [0.1cm]

\hline
\end{tabular}
\label{table_Eaton}
}
\end{table}

\begin{table}[H]
\caption{Data for the smoothed cube in 
three dimensions with timings in seconds.}

\vspace{.2in}

{\small
\centering
\begin{tabular}{cccccccc}  \\ [-0.3cm]
\hline
Size $(\lambda_{0})$ & $N_{tot}$ & $N$ & $E_2$ & $E_{inf}$ & $N_{matvec}$ & $T_{solve}$ & $T_{precomp}$ \\
\hline \\ [-0.3cm]
$1$ & $125000$ & $50$ & $4.08\times10^{-8}$ & $6.09\times10^{-8}$ & $15$ & $1.17$ & $1.07$ \\ [0.1cm]
$1$ & $343000$ & $70$ & $1.01\times10^{-10}$ & $1.25\times10^{-10}$ & $15$ & $2.05$ &  $2.23$\\ [0.1cm]
$1$ & $1000000$ & $100$ & $6.4\times10^{-14}$ & $7.91\times10^{-14}$ & $15$ & $5.04$ &  $6.44$\\ [0.1cm]
\hline
$20$ & $343000$ & $70$ & $7.23\times10^{-4}$ & $9.23\times10^{-4}$ & $449$ & $44.2$ & $2.11$\\ [0.1cm]
$20$ & $1000000$ & $100$ & $3.84\times10^{-8}$ & $5.22\times10^{-8}$ & $441$ & $124$ & $6.15$ \\ [0.1cm]
$20$ & $3375000$ & $150$ & $8.57\times10^{-13}$ & $1.64\times10^{-13}$ & $434$ & $411$ &  $20.1$\\ [0.1cm]
\hline
$40$ & $3375000$ & $150$ & $4.6\times10^{-7}$ & $4.34\times10^{-7}$ & $895$ & $891$ & $19.3$\\ [0.1cm]
$40$ & $8000000$ & $200$ & $6.76\times10^{-12}$ & $1.15\times10^{-11}$ & $907$ & $1957$ & $43.9$ \\ [0.1cm]
$40$ & $15625000$ & $250$ & $3.38\times10^{-12}$ & $5.76\times10^{-12}$ & $905$ & $3428$ &  $99.1$\\ [0.1cm]
\hline
$60$ & $3375000$ & $150$ & $4.76\times10^{-1}$ & $4.78\times10^{-1}$ & $8548$ & $8534$ & $19.8$\\ [0.1cm]
$60$ & $8000000$ & $200$ & $5.67\times10^{-6}$ & $4.99\times10^{-6}$ & $1471$ & $3538$ & $48.5$ \\ [0.1cm]
$60$ & $15625000$ & $250$ & $3.38\times10^{-9}$ & $5.86\times10^{-9}$ & $1505$ & $5688$ &  $97.8$\\ [0.1cm]
\hline
$80$ & $15625000$ & $250$ & $4.86\times10^{-5}$ & $4.74\times10^{-5}$ & $2988$ & $11332$ &  $97.5$\\ [0.1cm]
\hline

\end{tabular}
\label{table_SsquareFB}
}
\end{table}

\section{Conclusions}

We have presented a simple fast algorithm for computing volume potentials 
based on translation-invariant free-space Green's functions with 
compactly supported, smooth source distributions. By truncating the 
range of interaction in physical space, the Fourier transform of the kernel
becomes an entire function, so that the trapezoidal rule yields superalgebraic
convergence. Moreover, the transforms of the truncated kernels
can be computed analytically. The principal advantages of our approach are
that the standard FFT can be used for acceleration and that matrix 
entries for a high-order
accurate Nystr\"{o}m discretization are available ``on the fly". The latter
is important in constructing hierarchical direct solvers 
or incomplete LU preconditioners.

We have illustrated the performance of the scheme on a variety of problems
in two and three dimensions. For non-oscillatory problems, iterative methods
are quite effective at solving variable coefficient partial differential 
equations when recast as volume integral equations. For scattering problems,
it is well-known that the condition number grows with the size of the domain 
(measured in wavelengths). For problems up to approximately 
one hundred wavelengths in size, however, iterative schemes appear to be
viable without preconditioning.

We will explore the use of these methods for 
full electromagnetic scattering problems in three dimensions in future work.


\section*{Acknowledgments}
This work was supported in part by the
Applied Mathematical Sciences Program of the U.S. Department of Energy
under Contract DEFGO288ER25053 and
by the Office of the Assistant Secretary of Defense for Research and
Engineering and AFOSR under NSSEFF Program Award FA9550-10-1-0180.
The authors would like to thank 
Lise-Marie Imbert-Gerard and Carlos Borges for several useful
conversations.
The authors would also like to thank the ASIC 
(Area de Sistemas de Informaci\'on y Comunicaciones) and Francisco 
Jos\'e Rosich Viana for technical support and access to the UPV 
super-computing cluster RIGEL.

\bibliographystyle{unsrt}
\bibliography{referencias}

\begin{thebibliography}{10}

\bibitem{Ethridge2001}
F.~Ethridge and L.~Greengard.
\newblock {A New Fast-Multipole Accelerated Poisson Solver in Two Dimensions}.
\newblock {\em SIAM Journal on Scientific Computing}, 23(3):741--760, January
  2001.

\bibitem{Langston2008}
H.~Langston, L.~Greengard, and D.~Zorin.
\newblock A free-space adaptive fmm-based pde solver in three dimensions.
\newblock {\em Comm. Appl. Math. and Comp. Sci.}, 6:79--122, 2011.

\bibitem{Biros2}
D.~Malhotra and G.~Biros.
\newblock A parallel kernel independent fmm for particle and volume potentials.
\newblock {\em Communications in Computational Physics}, 18:808--830, 2015.

\bibitem{mcq}
P.~Mccorquodale, P.~Colella, G.~T. Balls, and S.~B. Baden.
\newblock A scalable parallel poisson solver in three dimensions with
  infinite-domain boundary conditions.
\newblock In {\em In 7th International Workshop on High Performance Scientific
  and Engineering Computing}, pages 814--822, 2005.

\bibitem{Aguilar2002}
J.~C. Aguilar and Y.~Chen.
\newblock {High-Order Corrected Trapezoidal Quadrature Rules for Functions with
  a Logarithmic Singularity in 2-D}.
\newblock {\em Computers {\&} Mathematics with Applications}, 44:1031--1039,
  2002.

\bibitem{Aguilar2005}
J.~C. Aguilar and Y.~Chen.
\newblock {High-order corrected trapezoidal quadrature rules for the coulomb
  potential in three dimensions}.
\newblock {\em Computers and Mathematics with Applications}, 49(4):625--631,
  2005.

\bibitem{Singular_regular_1}
J.~T. Beale and M.-C. Lai.
\newblock A method for computing nearly singular integrals.
\newblock {\em SIAM J. Numer. Anal}, 38:1902--25, 2001.

\bibitem{Duan2009}
R.~Duan and V.~Rokhlin.
\newblock {High-order quadratures for the solution of scattering problems in
  two dimensions}.
\newblock {\em Journal of Computational Physics}, 228(6):2152--2174, apr 2009.

\bibitem{goodman_1990}
J.~Goodman, T.~Y. Hou, and J.~Lowengrub.
\newblock The convergence of the point vortex method for the 2-d euler
  equations.
\newblock {\em Communications on Pure and Applied Mathematics}, 43:415--430,
  1990.

\bibitem{lowengrub_1993}
J.~Lowengrub, M.~Shelley, and B.~Merriman.
\newblock High-order and efficient methods for the vorticity formulation of the
  euler equations.
\newblock {\em {SIAM} Journal on Scientific Computing}, 14:1107--1142, 1993.

\bibitem{jiangbao}
S.~Jiang, L.~Greengard, and W.~Bao.
\newblock Fast and accurate evaluation of nonlocal coulomb and dipole-dipole
  interactions via the nonuniform fft.
\newblock {\em SIAM Journal on Scientific Computing}, 36:B777--B794, 2014.

\bibitem{Trefethen}
L.~N. Trefethen.
\newblock {\em Spectral methods in MATLAB}.
\newblock SIAM, Philadelphia, 2000.

\bibitem{Vainikko2000}
G.~Vainikko.
\newblock {Fast solvers of the Lippmann-Schwinger equation}.
\newblock {\em Direct and Inverse Problems of Mathematical Physics},
  5:423--440, 2000.

\bibitem{Hormander}
L.~H\"{o}rmander.
\newblock {\em Linear Partial Differential Equations}.
\newblock Springer, Berlin, 1976.

\bibitem{DelPhi}
L.~Li, C.~Li, Z.~Zhang, and E.~Alexov.
\newblock On the dielectric ÒconstantÓ of proteins: smooth dielectric function
  for macromolecular modeling and its implementation in delphi.
\newblock {\em Journal of chemical theory and computation}, 9(4):2126--2136,
  2013.

\bibitem{Nicholls}
J.~A. Grant, B.~T. Pickup, and A.~A. Nicholls.
\newblock {A smooth permittivity function for Poisson-Boltzmann solvation
  methods}.
\newblock {\em Journal of Computational Chemistry}, 230(22):608--640, 2001.

\bibitem{luneburg1964}
R.~K. Luneburg and M.~Herzberger.
\newblock {\em Mathematical theory of optics}.
\newblock Univ of California Press, 1964.

\bibitem{Eaton}
A.~J. Danner and U.~Leonhardt.
\newblock Lossless design of an {E}aton lens and invisible sphere by
  transformation optics with no bandwidth limitation.
\newblock In {\em Conference on Lasers and Electro-Optics}, page JThC4. Optical
  Society of America, 2009.

\end{thebibliography}

\end{document}